\documentclass[preprint,12pt]{elsarticle}

\usepackage{iftex}
\ifPDFTeX
  \usepackage[utf8]{inputenc}
\fi
\usepackage[T1]{fontenc}

\usepackage{BOONDOX-cal}
\usepackage{amssymb}
\usepackage{amsthm}
\usepackage{latexsym}
\usepackage{amsmath}
\usepackage{color}
\usepackage{graphicx}
\usepackage{indentfirst}
\usepackage{mathrsfs}
\usepackage{bm}
\usepackage[numbers,sort&compress]{natbib}
\usepackage{hyperref}
\usepackage{subcaption}
\usepackage{float}     
\usepackage{placeins}
\usepackage[titletoc]{appendix}

\hypersetup{
  colorlinks=true,
  linkcolor=blue,
  anchorcolor=blue,
  citecolor=blue,
  urlcolor=blue,
  pdftitle={Title of the Article},
  pdfauthor={Author One; Author Two}
}

\numberwithin{equation}{section}
\allowdisplaybreaks

\newtheorem{theorem}{\color{black}\indent Theorem}[section]

\newtheorem{assumption}{\color{black}\indent Assumption}[section]
\newtheorem{lemma}{\color{black}\indent Lemma}[section]
\newtheorem{proposition}{\color{black}\indent Proposition}[section]
\newtheorem{definition}{\color{black}\indent Definition}[section]
\newtheorem{remark}{\color{black}\indent Remark}[section]
\newtheorem{corollary}{\color{black}\indent Corollary}[section]

\newtheorem*{lemma*}{Lemma}

\newcommand{\TT}{\mathbb{T}}
\providecommand{\T}{\mathbb{T}} 

\newcommand{\RR}{\mathbb{R}}

\newcommand{\dd}{\,\mathrm{d}}

\newcommand{\avgtheta}[1]{\left\langle #1 \right\rangle_{\theta}}
\newcommand{\abs}[1]{\left|#1\right|}

\newcommand{\1}{\mathbf{1}}

\newcommand{\grad}{\nabla}
\DeclareMathOperator{\Div}{div}
\newcommand{\norm}[1]{\left\|#1\right\|}
\providecommand{\R}{\mathbb{R}}
\providecommand{\Z}{\mathbb{Z}}

\journal{xxxxxxxx}

\begin{document}

\biboptions{numbers,sort&compress}

\begin{frontmatter}

\title{\textbf{Statistical Ensemble Deviation Estimates for Nearly Integrable Hamiltonian Systems}}

\author[author1]{Xinyu Liu}
\ead{liuxy595@jlu.edu.cn}

\author[author1,author2]{Yong Li\corref{cor1}}
\ead{liyong@jlu.edu.cn}

\address[author1]{College of Mathematics, Jilin University, Changchun 130012, P. R. China.}
\address[author2]{Center for Mathematics and Interdisciplinary Sciences, Northeast Normal University, Changchun 130024, P. R. China.}

\cortext[cor1]{Corresponding author.}

\begin{abstract}
This paper studies quantitative deviation bounds for statistical ensembles evolving under the one-parameter flow of a nearly integrable Hamiltonian system. Combining Nekhoroshev-type stability estimates with phase-mixing arguments, we obtain, for any observable $G$, an explicit upper bound on the deviation of the ensemble average $\langle G\rangle_t$ from its angular average $\langle \avgtheta{G}\rangle_{0}$ over exponentially long time scales. The bound separates contributions from the resonant neighborhood via a probability-mass term, and from the nonresonant region via a traceable $1/t$ mixing constant $C_G$, a high-frequency Fourier tail, and an explicit normal-form remainder error.
\end{abstract}


\end{frontmatter}


\section{Introduction}\label{sec:introduction}

The evolution of statistical ensembles under Hamiltonian flows is a basic theme in classical statistical mechanics and nonlinear dynamics. 
Given an initial probability density $f_0(\theta,I)$ on phase space and an observable $G(\theta,I)$, the ensemble expectation at time $t$ is
\[
\langle G\rangle_t:=\int_{\TT^n\times B} G\big(\Phi_\varepsilon^t(\theta,I)\big)\,f_0(\theta,I)\,\dd\theta\,\dd I,
\]
where $\Phi_\varepsilon^t$ denotes the one-parameter Hamiltonian flow.
A central question is whether and how $\langle G\rangle_t$ approaches an equilibrium value as $t$ grows, and how to quantify the deviation from such an equilibrium.

In the integrable setting, action--angle coordinates yield the explicit dynamics
$I(t)\equiv I(0)$ and $\theta(t)=\theta(0)+t\omega(I)$, so that the long-time behavior of ensembles is governed by phase mixing induced by the frequency map $\omega$.
Under suitable regularity and nondegeneracy assumptions, Mitchell~\cite{Mitchell} established a weak convergence to equilibrium: as $t\to\infty$,
\[
\langle G\rangle_t \longrightarrow \langle \avgtheta{G}\rangle_0,
\qquad 
\avgtheta{G}(I):=\frac{1}{(2\pi)^n}\int_{\TT^n}G(\theta,I)\,\dd\theta.
\]

Beyond the fully integrable setting, it is natural to ask to what extent phase-mixing-based ensemble behavior persists under perturbations.
Several intermediate classes of perturbed integrable models have been studied in this direction.
For instance, we considered integrable Hamiltonian systems with transition terms~\cite{LiuLi2025AlmostPeriodic}, with external forcing~\cite{Liu2024PeriodicForced}, and with stochastic perturbations~\cite{LiuZhangLi2026LLNCLT}, and derived quantitative estimates for ensemble expectations under suitable regularity assumptions.
Moreover, besides convergence of ensemble means, one may also ask about fluctuations; central-limit-theorem-type results for appropriately normalized averages have been proved in related settings~\cite{LiuZhangLi2026LLNCLT}.
The present paper addresses a particularly important and more delicate regime: deterministic nearly integrable Hamiltonians, where perturbations destroy the exact torus translation structure and create a resonant/nonresonant phase-space geometry, so that global long-time behavior must be quantified through stability theory rather than pure integrable mixing.

In this work we consider a nearly integrable Hamiltonian system of the form
\begin{equation*}
H_\varepsilon(\theta,I)=h(I)+f_\varepsilon(\theta,I),\qquad 0<\varepsilon\ll 1,
\end{equation*}
where the perturbation couples angles and actions and creates a highly non-uniform phase-space structure, especially near low-order resonances.
Consequently, one should not expect a global convergence statement analogous to the integrable case for arbitrary initial ensembles.
A more robust goal is to obtain explicit deviation estimates for
\begin{equation*}
\big|\langle G\rangle_t-\langle \avgtheta{G}\rangle_0\big|
\end{equation*}
over the long (typically exponentially long) time scales that are relevant to nearly integrable stability.

Two classical toolboxes naturally come into play for nearly integrable Hamiltonian systems: KAM theory and Nekhoroshev theory.

KAM theory shows that, for sufficiently small perturbations and under Diophantine-type nonresonance conditions, a large (typically Cantor-like) family of invariant Lagrangian tori persists, and the dynamics on each such torus is conjugate to a linear quasi-periodic flow \cite{CongKupperLiYou2000,BertiKappelerMontalto2021,BertiHassainiaMasmoudi2023,BertiFranzoiMaspero2024,HassainiaHmidiRoulley2024}.
This provides a refined description of long-time quasi-periodic motions, but it is not, by itself, designed to produce a quantitative \emph{relaxation/deviation bound} for ensemble expectations with a general initial density.
Indeed, on individual invariant tori there is no genuine mixing; any decay of ensemble correlations must come from averaging over families of tori (i.e., phase mixing across actions).
To turn the KAM picture into an explicit deviation estimate, one would need in addition (i) quantitative control of the probability mass carried by the complement of the KAM set (weighted by the given $f_0$), (ii) uniform bounds on the torus-wise conjugacies (and the induced distortion on densities/observables), and (iii) mixing estimates whose constants remain stable across the surviving Cantor family.
These requirements go beyond the standard qualitative KAM output and would require additional bookkeeping tailored to the ensemble.

By contrast, Nekhoroshev theory provides \emph{uniform} (in the initial condition) stability estimates over very long time intervals, and this uniformity is particularly suitable for ensemble averages.
Nekhoroshev theory also remains an active area with continuing refinements on time scales and deviation bounds \cite{ZhangZhang2017,Poschel1993,YuLi2025NekhoroshevRandom,BambusiFeolaMontalto2024,FeolaMassetti2023,Barbieri2025SemiAlgebraic}.
Although Nekhoroshev stability is finite-time in nature, the resulting exponentially long windows are already sufficient for many applications, and---crucially---can be made quantitative in a way compatible with statistical ensembles.

The main idea of this paper is to combine Nekhoroshev-type stability with quantitative phase-mixing estimates.
We partition the action domain into a neighborhood $N(\varepsilon)$ of low-order resonances and its complement $D(\varepsilon)$.
The contribution of $N(\varepsilon)$ is controlled by the initial probability mass
\[
\mathsf P_{\mathrm{res}}(\varepsilon):=\int_{N(\varepsilon)\times\TT^n} f_0(\theta,I)\,\dd\theta\,\dd I.
\]
On the nonresonant region $D(\varepsilon)$ we apply an analytic normal-form transformation, which conjugates the original system to a Hamiltonian of the form
$\widetilde H_\varepsilon=h_\varepsilon(I)+r_\varepsilon(\theta,I)$, where $h_\varepsilon$ is integrable and the remainder $r_\varepsilon$ is exponentially small in the normal-form order.
We then establish a ready-to-use quantitative phase-mixing estimate for the integrable principal flow (with an explicit $1/t$ rate and a traceable constant), and control the contribution of $r_\varepsilon$ as an error term.

As a result, for any $t\neq 0$ we obtain an explicit deviation bound of the form
\begin{align*}
  \begin{aligned}
\big|\langle G\rangle_t-\langle \avgtheta{G}\rangle_0\big|
\ &\le\
2\|G\|_\infty\,\mathsf P_{\mathrm{res}}(\varepsilon)
+\frac{C_G(K(\varepsilon);D(\varepsilon))}{|t|}\\
&+\mathcal R_{>K(\varepsilon)}
+\mathcal E_{\mathrm{nf}}(t,\varepsilon)
+\mathcal E_{\mathrm{eq}}(\varepsilon),
  \end{aligned}
\end{align*}
where the terms on the right-hand side separately quantify the resonant probability mass, the nonresonant phase mixing, the high-frequency Fourier tail, and the normal-form and coordinate-change errors.
In particular, when $\|r_\varepsilon\|_\infty\lesssim \varepsilon e^{-cK(\varepsilon)}$, the normal-form error remains exponentially small on exponentially long time windows $|t|\le e^{\sigma K(\varepsilon)}$ with $0<\sigma<c/2$, consistent with Nekhoroshev-type stability scales.

The remainder of the paper is organized as follows.
Section~\ref{sec:preliminaries} introduces the setting, the resonant/nonresonant decomposition, and the basic ensemble notation.
Section~\ref{sec:main-results} first derives quantitative phase-mixing estimates for integrable systems and then embeds them into the nearly integrable framework via normal forms, leading to the main deviation theorem and its corollaries.
The final section summarizes the results and discusses perspectives.

\section{Preliminaries}\label{sec:preliminaries}
In this section, we introduce the system under consideration and collect the basic notation and definitions used throughout the paper.

Consider a nearly integrable Hamiltonian system
\begin{equation}\label{EQ:NIH}
  H_\varepsilon(\theta,I)=h(I)+f_{\varepsilon}(\theta,I),
\end{equation}
in which $(\theta,I)\in \TT^n\times B(0,R)\subset \TT^n\times\RR^n,$ $h$ denotes the integrable part, and $f_{\varepsilon}$ is a small perturbation.
Denote by $\Phi_\varepsilon^t$ the Hamiltonian flow generated by $H_\varepsilon$.

Let the frequency map be defined by
$$\omega(I):=\nabla h(I).$$
For each $k\in\mathbb Z^n\setminus\{0\}$, define the resonance hyperplane
$$R_k:=\{\omega\in\mathbb R^n:\ k\cdot \omega=0\}.$$
\begin{definition}
  Given parameters $K(\varepsilon)>0$ and $\alpha(\varepsilon)>0$, define
\[
N(\varepsilon):=\Big\{I\in B(0,R): d\Big(\omega(I),\bigcup_{0<\abs{k}\le K(\varepsilon)}R_k\Big)<\alpha(\varepsilon)\Big\},
\quad
D(\varepsilon):=B(0,R)\setminus N(\varepsilon).
\]
\end{definition}

\begin{remark}
In \cite{ZhangZhang2017}, to construct a global stability partition, the authors fix $0<\beta<1$, set $L=12s_0$, and define
\[
K(\varepsilon)=-L\log\varepsilon,\qquad r(\varepsilon)=\beta^{-1}\varepsilon^{1/2},\qquad
\alpha(\varepsilon)=\beta^{-1}r(\varepsilon)K(\varepsilon).
\]
The results in this paper can be stated for general choices of $(K(\varepsilon),\alpha(\varepsilon))$.
\end{remark}

\begin{definition}
Let $f_0(\theta,I)\ge 0$ be an initial probability density on $\TT^n\times B$, normalized by
\[
\int_{\TT^n\times B} f_0\,\dd\theta\,\dd I=1.
\]
For an observable $G(\theta,I)$, define the ensemble expectation at time $t$ by
\[
\langle G\rangle_t:=\int_{\TT^n\times B} G\big(\Phi_\varepsilon^t(\theta_0,I_0)\big)\,f_0(\theta_0,I_0)\,\dd\theta_0\,\dd I_0.
\]
Define the angular average
\[
\avgtheta{G}(I):=\frac{1}{(2\pi)^n}\int_{\TT^n}G(\theta,I)\,\dd\theta,
\]
and its initial expectation
\[
\langle \avgtheta{G}\rangle_0:=\int_{\TT^n\times B}\avgtheta{G}(I)\,f_0(\theta,I)\,\dd\theta\,\dd I.
\]
\end{definition}

\begin{definition}
Let the initial marginal density be
\[
\rho_0(I):=\frac{1}{(2\pi)^n}\int_{\TT^n} f_0(\theta,I)\,\dd\theta.
\]
Define the initial probability mass of the resonant region by
\[
\mathsf{P}_{\mathrm{res}}(\varepsilon):=\int_{N(\varepsilon)\times \TT^n} f_0(\theta,I)\,\dd\theta\,\dd I
=\int_{N(\varepsilon)} \rho_0(I)\,\dd I.
\]
\end{definition}
\begin{remark}
This term quantifies how much probability mass of the initial ensemble lies in a neighborhood of low-order resonances; in the following paper it appears in the form $2\|G\|_{\infty}\,\mathsf{P}_{\mathrm{res}}(\varepsilon)$.
\end{remark}
\section{Main results}\label{sec:main-results}
This section consists of two parts. First, building on results of C.~Mitchell \cite{Mitchell}, we establish quantitative estimates for the relevant convergence results. Second, we embed the integrable mixing estimate into the nearly integrable setting, thereby deriving the main results of the paper.
\subsection{Lemmas and quantitative estimates for integrable systems}
First, we present two key lemmas. The first provides a quantitative version of the Riemann--Lebesgue lemma (with an explicit expression for $u$); the second bounds $\|u\|_{L^{1}}$ in terms of $\gamma:=\inf |\nabla\phi|$, $\|\nabla^{2}\phi\|_{\infty}$, and the $W^{1,1}$-norm of $a$. More precisely, we have the following.

\begin{lemma}\label{lem:MitchellIBP}
Let $\Omega\subset\RR^n$ be an open set, $a\in C_c^1(\Omega)$, and let $\phi\in C^2(\Omega)$ be real-valued with $\grad\phi\neq 0$ on $\Omega$.
For any $\lambda\in\R\setminus\{0\}$, define
\[
I(\lambda):=\int_\Omega a(x)\,e^{i\lambda \phi(x)}\,\dd x.
\]
Then the following identity holds:
\[
I(\lambda)=-\frac{1}{i\lambda}\int_\Omega u(x)\,e^{i\lambda\phi(x)}\,\dd x,
\qquad
u(x):=\Div\!\left(a(x)\frac{\grad\phi(x)}{\abs{\grad\phi(x)}^2}\right),
\]
and consequently,
\[
\abs{I(\lambda)}\le \frac{1}{\abs{\lambda}}\norm{u}_{L^1(\Omega)}.
\]
\end{lemma}

\begin{proof}
By the chain rule,
\[
\grad e^{i\lambda\phi}= i\lambda(\grad\phi)\,e^{i\lambda\phi}.
\]
Hence, on $\Omega$,
\[
e^{i\lambda\phi}=\frac{1}{i\lambda}\frac{\grad\phi}{\abs{\grad\phi}^2}\cdot \grad e^{i\lambda\phi}.
\]
Substituting this into $I(\lambda)$ yields
\[
I(\lambda)=\frac{1}{i\lambda}\int_\Omega a\,\frac{\grad\phi}{\abs{\grad\phi}^2}\cdot \grad e^{i\lambda\phi}\,\dd x.
\]
Since $a\in C_c^1(\Omega)$ has compact support (hence vanishes on $\partial\Omega$), we may integrate by parts and drop the boundary term:
\[
I(\lambda)=-\frac{1}{i\lambda}\int_\Omega
\Div\!\left(a\frac{\grad\phi}{\abs{\grad\phi}^2}\right)e^{i\lambda\phi}\,\dd x
=-\frac{1}{i\lambda}\int_\Omega u\,e^{i\lambda\phi}\,\dd x.
\]
Taking absolute values and using the triangle inequality gives
$\abs{I(\lambda)}\le \abs{\lambda}^{-1}\norm{u}_{L^1}$.
\end{proof}

\begin{lemma}\label{lem:uL1bound}
Under the assumptions of Lemma~\ref{lem:MitchellIBP}, assume in addition that
\[
\gamma:=\inf_{x\in\Omega}\abs{\grad\phi(x)}>0,\qquad
M:=\norm{\grad^2\phi}_{L^\infty(\Omega)}<\infty.
\]
Then
\[
\norm{u}_{L^1(\Omega)}
\le
\frac{1}{\gamma}\,\norm{\grad a}_{L^1(\Omega)}
+\frac{(n+2)M}{\gamma^2}\,\norm{a}_{L^1(\Omega)}.
\]
\end{lemma}

\begin{proof}
By the product rule,
\[
u=\Div\!\left(a\frac{\grad\phi}{\abs{\grad\phi}^2}\right)
=
\grad a\cdot\frac{\grad\phi}{\abs{\grad\phi}^2}
+
a\,\Div\!\left(\frac{\grad\phi}{\abs{\grad\phi}^2}\right).
\]
Hence,
\[
\abs{u}\le
\abs{\grad a}\frac{1}{\abs{\grad\phi}}
+
\abs{a}\,\abs{\Div(\grad\phi/\abs{\grad\phi}^2)}.
\]
For the first term, using $\abs{\grad\phi}\ge\gamma$ we have
$\abs{\grad a}/\abs{\grad\phi}\le \abs{\grad a}/\gamma$.
For the second term, note the identity
\[
\Div\!\left(\frac{\grad\phi}{\abs{\grad\phi}^2}\right)
=
\frac{\Delta\phi}{\abs{\grad\phi}^2}
-2\frac{(\grad\phi)^\top(\grad^2\phi)\grad\phi}{\abs{\grad\phi}^4}.
\]
Since $\abs{\Delta\phi}\le n\norm{\grad^2\phi}$ and
$\abs{(\grad\phi)^\top(\grad^2\phi)\grad\phi}\le \norm{\grad^2\phi}\abs{\grad\phi}^2$,
we obtain the pointwise bound
\[
\abs{\Div(\grad\phi/\abs{\grad\phi}^2)}
\le \frac{(n+2)\norm{\grad^2\phi}}{\abs{\grad\phi}^2}
\le \frac{(n+2)M}{\gamma^2}.
\]
Combining the two estimates and integrating over $\Omega$ yields the desired inequality.
\end{proof}
Next, we work within the \emph{integrable} action--angle framework. Our goal is to provide a ``ready-to-use'' $1/t$ phase-mixing estimate (with a traceable constant $C_G$), which will later be applied to the integrable normal-form principal part in the nearly integrable setting.

Consider
\begin{equation}\label{ass:int}
\dot I=0,\qquad \dot\theta=\omega(I),\qquad (I,\theta)\in \Omega\times\TT^n,
\end{equation}
whose flow is given by $\Psi_t(I,\theta)=(I,\theta+t\omega(I))$.
Let the initial density satisfy $f_0\in C_c^1(\Omega\times\T^n)$, and let the observable be $G\in C(\Omega\times\T^n)$.

\begin{definition}\label{def:fourier}
Expand in Fourier series with respect to $\theta\in\T^n$:
\[
G(\theta,I)=\sum_{k\in\Z^n}G_k(I)e^{ik\cdot\theta},\qquad
f_0(\theta,I)=\sum_{k\in\Z^n}f_{0,k}(I)e^{ik\cdot\theta}.
\]
Given $K\ge 1$, define the high-frequency tail term
\[
\mathcal R_{>K}(G,f_0;\Omega):=(2\pi)^n\sum_{\abs{k}>K}\int_\Omega \abs{G_k(I)f_{0,-k}(I)}\,\dd I.
\]
\end{definition}

\begin{proposition}\label{prop:mixing}
For the system~\eqref{ass:int}, fix $K\ge 1$ and set
\[
a_k(I):=G_k(I)\,f_{0,-k}(I),\qquad
\phi_k(I):=k\cdot\omega(I).
\]
Assume that for each $0<\abs{k}\le K$, one has $a_k\in W^{1,1}(\Omega)$ and
\[
\gamma_k(\Omega):=\inf_{I\in\Omega}\abs{\grad_I\phi_k(I)}>0.
\]
Then for any $t\neq 0$,
\[
\abs{\langle G\rangle_t-\langle \avgtheta{G}\rangle_0}
\le
\frac{C_G(K;\Omega)}{\abs{t}}+\mathcal R_{>K}(G,f_0;\Omega),
\]
where
\[
C_G(K;\Omega):=(2\pi)^n\sum_{0<\abs{k}\le K}\norm{u_k}_{L^1(\Omega)},
\qquad
u_k(I):=\Div_I\!\left(a_k(I)\frac{\grad_I\phi_k(I)}{\abs{\grad_I\phi_k(I)}^2}\right).
\]

Moreover, if $\omega\in C^2(\Omega)$ and we define
\[
M_k(\Omega):=\norm{\grad_I^2\phi_k}_{L^\infty(\Omega)}
=\norm{\grad_I^2(k\cdot\omega)}_{L^\infty(\Omega)},
\]
then Lemma~\ref{lem:uL1bound} yields the explicit traceable bound
\[
C_G(K;\Omega)
\le
(2\pi)^n\sum_{0<\abs{k}\le K}
\left[
\frac{\norm{\grad_I a_k}_{L^1(\Omega)}}{\gamma_k(\Omega)}
+
\frac{(n+2)M_k(\Omega)}{\gamma_k(\Omega)^2}\norm{a_k}_{L^1(\Omega)}
\right].
\]
\end{proposition}

\begin{proof}
For the integrable flow $\Psi_t(I,\theta)=(I,\theta+t\omega(I))$, we have
\[
\langle G\rangle_t
=\int_{\Omega\times\TT^n} G(I,\theta+t\omega(I))\,f_0(I,\theta)\,\dd\theta\,\dd I.
\]
Substituting the Fourier expansions and integrating in $\theta$ (retaining only the terms with $k+m=0$) gives
\[
\langle G\rangle_t
=(2\pi)^n\sum_{k\in\Z^n}\int_\Omega a_k(I)\,e^{it\phi_k(I)}\,\dd I.
\]
The term $k=0$ is precisely $\langle \avgtheta{G}\rangle_0$.
Hence the deviation satisfies
\[
\langle G\rangle_t-\langle \avgtheta{G}\rangle_0
=(2\pi)^n\sum_{k\ne 0}\int_\Omega a_k(I)\,e^{it\phi_k(I)}\,\dd I.
\]
Split the sum into the low-frequency part $\abs{k}\le K$ and the high-frequency part $\abs{k}>K$.
The high-frequency contribution is estimated in absolute value, yielding $\mathcal R_{>K}$.
For each $0<\abs{k}\le K$, apply Lemma~\ref{lem:MitchellIBP} to the oscillatory integral
(with $\lambda=t$, $a=a_k$, and $\phi=\phi_k$) to obtain
\[
\abs{\int_\Omega a_k(I)e^{it\phi_k(I)}\,\dd I}\le \frac{\norm{u_k}_{L^1(\Omega)}}{\abs{t}}.
\]
Summing over $0<\abs{k}\le K$ yields the stated estimate, and the explicit bound on $C_G$ follows from Lemma~\ref{lem:uL1bound}.
\end{proof}

\subsection{Deviation of statistical ensembles of nearly integrable Hamiltonian}

Now we can present the main results of this paper. 
The core idea is as follows: First, we decompose the initial ensemble according to the partition $D(\varepsilon)\cup N(\varepsilon)$. The contribution from the resonant neighborhood $N(\varepsilon)$ is then bounded by the probability-mass term $2\|G\|_{\infty}\,\mathsf P_{\rm res}(\varepsilon)$, and on $D(\varepsilon)$, we use a normal-form transformation to rewrite the nearly integrable system as an integrable principal part plus a small remainder. We then apply Proposition~\ref{prop:mixing} to the integrable principal part and control the contribution of the remainder as an error term.

To avoid notational ambiguity, we use the following conventions. We write $\Phi_\varepsilon^t$ for the flow generated by the original Hamiltonian $H_\varepsilon$, $\Phi_\varepsilon^{\rm nf}$ for the normal-form symplectic transformation, and $\widetilde H_\varepsilon:=H_\varepsilon\circ\Phi_\varepsilon^{\rm nf}=h_\varepsilon+r_\varepsilon$ for the normal-form Hamiltonian, whose flow is denoted by $\widetilde\Phi_\varepsilon^t$.

Next, we first state the pullback lemma for ensemble expectations in normal-form coordinates.

\begin{lemma}\label{lem:nf-pullback}
Let the phase space $M:=\T^n\times D(\varepsilon)$ be endowed with the standard symplectic form
\[
\omega_0:=\sum_{j=1}^n \dd I_j\wedge \dd\theta_j,
\]
and let the Liouville volume measure be
\[
\mu:=\frac{\omega_0^n}{n!}=\dd\theta\,\dd I.
\]
Assume that $H_\varepsilon$ is $C^1$ on an open neighborhood of $M$, and that its Hamiltonian vector field $X_{H_\varepsilon}$ generates a flow $\Phi_\varepsilon^t$ on $M$ (i.e., for all initial data in $M$ the flow exists on the time interval under consideration).

Suppose that there exists a $C^1$ symplectic diffeomorphism (symplectic transformation)
\[
\Phi_\varepsilon^{\rm nf}:M\to M,
\qquad (\Phi_\varepsilon^{\rm nf})^*\omega_0=\omega_0,
\]
and define the normal-form Hamiltonian
\[
\widetilde H_\varepsilon:=H_\varepsilon\circ \Phi_\varepsilon^{\rm nf}.
\]
Denote by $\widetilde\Phi_\varepsilon^t$ the flow generated by $\widetilde H_\varepsilon$ (assumed to exist on the same time interval).

For any measurable function $G$ and any initial density $f_0^D\in L^1(M,\mu)$ (e.g., $f_0^D=f_0\,\1_{D(\varepsilon)}$), set
\[
\widetilde G:=G\circ \Phi_\varepsilon^{\rm nf},\qquad
\widetilde f_0:=f_0^D\circ \Phi_\varepsilon^{\rm nf}.
\]
Then for any admissible time $t$,
\[
\int_{M} G(\Phi_\varepsilon^t(z))\,f_0^D(z)\,\dd\mu(z)
=
\int_{M} \widetilde G\big(\widetilde\Phi_\varepsilon^t(z)\big)\,
\widetilde f_0(z)\,\dd\mu(z).
\]
\end{lemma}
\begin{proof}
Since $\Phi_\varepsilon^{\rm nf}$ is a symplectic diffeomorphism, it satisfies $(\Phi_\varepsilon^{\rm nf})^*\omega_0=\omega_0$.
Hence, for the $n$-fold exterior power,
\begin{equation*}
(\Phi_\varepsilon^{\rm nf})^*(\omega_0^n)
=\big((\Phi_\varepsilon^{\rm nf})^*\omega_0\big)^n
=\omega_0^n,
\end{equation*}
and therefore $(\Phi_\varepsilon^{\rm nf})^*\mu=\mu$, i.e., $\Phi_\varepsilon^{\rm nf}$ preserves the measure $\mu$.
In particular, in standard coordinates this is equivalent to
\begin{equation*}
\det D\Phi_\varepsilon^{\rm nf}(z)=1\quad \text{a.e. } z\in M,
\end{equation*}
so the change-of-variables formula can be written as
\begin{equation*}
\int_M F(z)\,\dd\mu(z)
=\int_M F(\Phi_\varepsilon^{\rm nf}(\tilde z))\,\dd\mu(\tilde z)
\qquad(\forall\,F\in L^1).
\end{equation*}

Let $\Phi:=\Phi_\varepsilon^{\rm nf}$. We prove that
\[
\Phi_\varepsilon^t\circ \Phi=\Phi\circ \widetilde\Phi_\varepsilon^t.
\tag{$\ast$}
\]
To this end, we first establish the transformation rule for Hamiltonian vector fields:
\[
D\Phi(z)\,X_{\widetilde H_\varepsilon}(z)=X_{H_\varepsilon}(\Phi(z))\qquad(\forall z\in M).
\tag{$\dagger$}
\]
Recall the definition of the Hamiltonian vector field: $i_{X_{H_\varepsilon}}\omega_0=\dd H_\varepsilon$ and
$i_{X_{\widetilde H_\varepsilon}}\omega_0=\dd \widetilde H_\varepsilon$.
Fix $z\in M$ and an arbitrary tangent vector $Y\in T_zM$. Using that $\Phi$ preserves the symplectic form, we obtain
\[
\omega_0\big(D\Phi(z)X_{\widetilde H_\varepsilon}(z),\,D\Phi(z)Y\big)
=
\omega_0\big(X_{\widetilde H_\varepsilon}(z),\,Y\big)
=
\dd\widetilde H_\varepsilon(z)[Y].
\]
On the other hand, since $\widetilde H_\varepsilon=H_\varepsilon\circ\Phi$, we have
\[
\dd\widetilde H_\varepsilon(z)[Y]
=
\dd(H_\varepsilon\circ\Phi)(z)[Y]
=
\dd H_\varepsilon(\Phi(z))[D\Phi(z)Y]
=
\omega_0\big(X_{H_\varepsilon}(\Phi(z)),\,D\Phi(z)Y\big).
\]
Therefore, for all $Y$,
\[
\omega_0\big(D\Phi(z)X_{\widetilde H_\varepsilon}(z)-X_{H_\varepsilon}(\Phi(z)),\,D\Phi(z)Y\big)=0.
\]
Since $D\Phi(z)$ is invertible and $\omega_0$ is nondegenerate, this implies
$D\Phi(z)X_{\widetilde H_\varepsilon}(z)=X_{H_\varepsilon}(\Phi(z))$, which is $(\dagger)$.

Now fix any initial condition $z_0\in M$, and set $\tilde z(t):=\widetilde\Phi_\varepsilon^t(z_0)$.
By definition, $\dot{\tilde z}(t)=X_{\widetilde H_\varepsilon}(\tilde z(t))$.
Define $z(t):=\Phi(\tilde z(t))$. Then, by the chain rule and $(\dagger)$,
\[
\dot z(t)=D\Phi(\tilde z(t))\,\dot{\tilde z}(t)
=D\Phi(\tilde z(t))X_{\widetilde H_\varepsilon}(\tilde z(t))
=X_{H_\varepsilon}(\Phi(\tilde z(t)))=X_{H_\varepsilon}(z(t)).
\]
Thus $z(t)$ solves the Hamiltonian system generated by $H_\varepsilon$, with $z(0)=\Phi(z_0)$.
By uniqueness of solutions to ordinary differential equations,
\[
z(t)=\Phi_\varepsilon^t(\Phi(z_0)).
\]
On the other hand, $z(t)=\Phi(\tilde z(t))=\Phi(\widetilde\Phi_\varepsilon^t(z_0))$.
Hence, for every $z_0$,
\[
\Phi_\varepsilon^t(\Phi(z_0))=\Phi(\widetilde\Phi_\varepsilon^t(z_0)),
\]
which yields the conjugacy relation $(\ast)$.

Let $M=\T^n\times D(\varepsilon)$. Starting from the left-hand side,
\[
\int_{M} G(\Phi_\varepsilon^t(z))\,f_0^D(z)\,\dd\mu(z).
\]
Perform the change of variables $z=\Phi(\tilde z)$. By the measure-preserving property from  $\dd\mu(z)=\dd\mu(\tilde z)$, we obtain
\[
\int_{M} G\big(\Phi_\varepsilon^t(\Phi(\tilde z))\big)\,f_0^D(\Phi(\tilde z))\,\dd\mu(\tilde z).
\]
By the conjugacy relation from $\Phi_\varepsilon^t\circ\Phi=\Phi\circ\widetilde\Phi_\varepsilon^t$, we have
\[
G\big(\Phi_\varepsilon^t(\Phi(\tilde z))\big)
=G\big(\Phi(\widetilde\Phi_\varepsilon^t(\tilde z))\big)
=(G\circ\Phi)\big(\widetilde\Phi_\varepsilon^t(\tilde z)\big)
=\widetilde G\big(\widetilde\Phi_\varepsilon^t(\tilde z)\big).
\]
Moreover,
\[
f_0^D(\Phi(\tilde z))=(f_0^D\circ\Phi)(\tilde z)=\widetilde f_0(\tilde z).
\]
Substituting these identities yields
\[
\int_{M} \widetilde G\big(\widetilde\Phi_\varepsilon^t(\tilde z)\big)\,\widetilde f_0(\tilde z)\,\dd\mu(\tilde z),
\]
which is exactly the right-hand side. This completes the proof.
\end{proof}

\begin{remark}\label{rem:domain}
In order for the change of variables to be performed over the same integration domain $M=\T^n\times D(\varepsilon)$, we assumed above that
$\Phi_\varepsilon^{\rm nf}:M\to M$ is a symplectic diffeomorphism.
In actual normal-form constructions, $\Phi_\varepsilon^{\rm nf}$ is typically defined only on a slightly larger domain and satisfies
$\Phi_\varepsilon^{\rm nf}(M_-)\subset M$ $($or $M\subset \Phi_\varepsilon^{\rm nf}(M_-)$$)$.
In this case, one may either take the integration domain to be $M_-$, or shrink $D(\varepsilon)$ appropriately so that the above self-mapping assumption holds.
Without such a restriction, the more general correct formulation is
\[
\int_{\Phi(M_-)} G(\Phi_\varepsilon^t(z))\,f_0(z)\,\dd\mu(z)
=
\int_{M_-} (G\circ \Phi_\varepsilon^{\rm nf})(\widetilde\Phi_\varepsilon^t(z))\,
(f_0\circ \Phi_\varepsilon^{\rm nf})(z)\,\dd\mu(z),
\]
where the integration domain on the left-hand side is taken to be $\Phi(M_-)$, while the right-hand side is integrated over $M_-$.
\end{remark}

\begin{lemma}\label{lem:nf-pullback1}
Let $f_0^D:=f_0\,\1_{D(\varepsilon)}$.
Assume that there exists a symplectic transformation $\Phi_\varepsilon^{\rm nf}$ on $D(\varepsilon)$, and set
\[
\widetilde G:=G\circ \Phi_\varepsilon^{\rm nf},\qquad
\widetilde f_0:=f_0^D\circ \Phi_\varepsilon^{\rm nf},\qquad
\widetilde H_\varepsilon:=H_\varepsilon\circ \Phi_\varepsilon^{\rm nf}.
\]
Then
\[
\langle G\rangle_t^D
:=\int_{\T^n\times D(\varepsilon)} G(\Phi_\varepsilon^t(z))\,f_0^D(z)\,\dd z
=
\int_{\T^n\times D(\varepsilon)} \widetilde G\big(\widetilde\Phi_\varepsilon^t(z)\big)\,
\widetilde f_0(z)\,\dd z,
\]
where $\widetilde\Phi_\varepsilon^t$ denotes the Hamiltonian flow generated by $\widetilde H_\varepsilon$.
\end{lemma}

\begin{proof}
The claim follows from the conjugacy relation
$\Phi_\varepsilon^t\circ \Phi_\varepsilon^{\rm nf}
=
\Phi_\varepsilon^{\rm nf}\circ \widetilde\Phi_\varepsilon^t$
together with the change of variables $z=\Phi_\varepsilon^{\rm nf}(\tilde z)$,
and the fact that a symplectic transformation preserves the Liouville volume (equivalently, its Jacobian determinant equals $1$).
\end{proof}

\begin{lemma}\label{lem:nf-error}
Let the normal-form Hamiltonian be
$\widetilde H_\varepsilon(\theta,I)=h_\varepsilon(I)+r_\varepsilon(\theta,I)$,
and denote by
\[
\Psi_\varepsilon^t(\theta,I):=(\theta+t\omega_\varepsilon(I),I),\qquad \omega_\varepsilon=\grad h_\varepsilon,
\]
the translation flow generated by the integrable principal part.
Assume that $\widetilde G\in C^1$, that $\widetilde f_0\ge 0$, and that $\int \widetilde f_0=1$.
Then there exists a constant $C_{\rm err}>0$ (depending on the geometry of $D(\varepsilon)$, a $C^1$ bound on $\omega_\varepsilon$, and the analyticity radius)
such that for any $t\in\R$,
\[
\abs{
\int \widetilde G(\widetilde\Phi_\varepsilon^t(z))\,\widetilde f_0(z)\,\dd z
-
\int \widetilde G(\Psi_\varepsilon^t(z))\,\widetilde f_0(z)\,\dd z
}
\le
C_{\rm err}\,\|\widetilde G\|_{C^1}\,(1+\abs{t}+t^2)\,\|r_\varepsilon\|_\infty.
\]
\end{lemma}

\begin{proof}
Let $z(t)=(\theta(t),I(t))$ be a solution of $\widetilde H_\varepsilon$, and let $z^0(t)=(\theta^0(t),I^0(t))$ be the solution of the integrable principal part,
where $I^0(t)\equiv I(0)$ and $\theta^0(t)=\theta(0)+t\omega_\varepsilon(I(0))$.
From the equations
\[
\dot I=-\partial_\theta r_\varepsilon,\qquad
\dot\theta=\omega_\varepsilon(I)+\partial_I r_\varepsilon
\]
we obtain
\[
\abs{I(t)-I(0)}\le \abs{t}\,\|\partial_\theta r_\varepsilon\|_\infty.
\]
Moreover,
\[
\theta(t)-\theta^0(t)
=
\int_0^t\big(\omega_\varepsilon(I(s))-\omega_\varepsilon(I(0))\big)\,\dd s
+
\int_0^t \partial_I r_\varepsilon(\theta(s),I(s))\,\dd s.
\]
Hence,
\begin{align*}
  \begin{aligned}
\abs{\theta(t)-\theta^0(t)}
&\le
\abs{t}\,\|\partial_I r_\varepsilon\|_\infty
+
\|\grad\omega_\varepsilon\|_\infty\int_0^{\abs{t}}\abs{I(s)-I(0)}\,\dd s\\
&\le
\abs{t}\,\|\partial_I r_\varepsilon\|_\infty
+\tfrac12 \abs{t}^2 \|\grad\omega_\varepsilon\|_\infty\,\|\partial_\theta r_\varepsilon\|_\infty.
  \end{aligned}
\end{align*}
If $r_\varepsilon$ is analytic on a fixed complex neighborhood, then Cauchy estimates yield
$\|\partial_\theta r_\varepsilon\|_\infty+\|\partial_I r_\varepsilon\|_\infty\le C\,\|r_\varepsilon\|_\infty$,
and therefore
\begin{equation*}
\abs{\theta(t)-\theta^0(t)}+\abs{I(t)-I^0(t)}\le C(1+\abs{t}+t^2)\|r_\varepsilon\|_\infty.
\end{equation*}
Finally, using
\begin{equation*}
\abs{\widetilde G(z(t))-\widetilde G(z^0(t))}
\le \|\widetilde G\|_{C^1}\big(\abs{\theta(t)-\theta^0(t)}+\abs{I(t)-I^0(t)}\big),
\end{equation*}
and integrating against $\widetilde f_0$ yields the desired estimate.
\end{proof}

\begin{assumption}\label{ass:nf}
Let $D(\varepsilon)$ be the nonresonant region defined in Section~2. Assume that there exists a real-analytic symplectic transformation
\[
\Phi_\varepsilon^{\rm nf}:\T^n\times D(\varepsilon)\to \T^n\times B(0,R)
\]
such that on $\T^n\times D(\varepsilon)$,
\[
H_\varepsilon\circ \Phi_\varepsilon^{\rm nf} = h_\varepsilon(I)+r_\varepsilon(\theta,I),
\]
where $h_\varepsilon$ depends only on $I$ (and is therefore integrable), and the remainder satisfies
\[
\norm{r_\varepsilon}_{\infty}\le C_{\mathrm{nf}}\varepsilon\,e^{-c_{\mathrm{nf}}K(\varepsilon)}.
\]
\end{assumption}

\begin{remark}
In the setting of Zhang--Zhang \cite{ZhangZhang2017}, the normal form on $D(\varepsilon)$ follows from a P\"oschel-type lemma (e.g., their Lemma~3.1):
the remainder satisfies $\|f_1\|\le e^{-Ks/6}\|f\|$. Hence one may take $c_{\mathrm{nf}}=s/6$ in the above bound.
When $K(\varepsilon)\sim \varepsilon^{-a}$, one obtains $\|r_\varepsilon\|\lesssim \varepsilon e^{-c\varepsilon^{-a}}$,
which matches the exponential time scale in Nekhoroshev theory.
\end{remark}

\begin{assumption}\label{ass:fourierReg}
Let $f_0^D:=f_0\,\1_{D(\varepsilon)}$, and define in normal-form coordinates
\[
\widetilde G:=G\circ \Phi_\varepsilon^{\rm nf},\qquad
\widetilde f_0:=f_0^D\circ \Phi_\varepsilon^{\rm nf}.
\]
Let $\widetilde G_k$ and $\widetilde f_{0,k}$ denote their Fourier coefficients.
Assume that for all $0<\abs{k}\le K(\varepsilon)$,
\begin{equation*}
  \begin{gathered}
\widetilde a_k(I):=\widetilde G_k(I)\,\widetilde f_{0,-k}(I)\in W^{1,1}\big(D(\varepsilon)\big),\\
\gamma_k\big(D(\varepsilon)\big):=\inf_{I\in D(\varepsilon)}\abs{\grad_I\big(k\cdot \omega_\varepsilon(I)\big)}>0,
  \end{gathered}
\end{equation*}
where $\omega_\varepsilon(I):=\grad h_\varepsilon(I)$ is the frequency map of the integrable principal part in normal form.
\end{assumption}

\begin{theorem}[Main theorem: deviation bound for nearly integrable ensembles under the $D(\varepsilon)$--$N(\varepsilon)$ decomposition]\label{thm:main}
Under Assumptions~\ref{ass:nf} and~\ref{ass:fourierReg}, denote the probability mass of the resonant region by
\[
\mathsf P_{\mathrm{res}}(\varepsilon)=\int_{N(\varepsilon)\times\T^n} f_0(\theta,I)\,\dd\theta\,\dd I.
\]
Then for any $t\neq 0$,
\begin{align*}
\begin{aligned}
\abs{\langle G\rangle_t-\langle \avgtheta{G}\rangle_0}
\ &\le\
2\norm{G}_{L^\infty}\,\mathsf P_{\mathrm{res}}(\varepsilon)
\ +\
\frac{C_{\widetilde G,\widetilde f_0}\big(K(\varepsilon);D(\varepsilon)\big)}{\abs{t}}
\ \\
&+\
\mathcal R_{>K(\varepsilon)}\big(\widetilde G,\widetilde f_0;D(\varepsilon)\big)
\ +\
\mathcal E_{\mathrm{nf}}(t,\varepsilon)
\ +\
\mathcal E_{\mathrm{eq}}(\varepsilon).
\end{aligned}
\end{align*}
Here:
\begin{itemize}
\item $C_{\widetilde G,\widetilde f_0}(K;D)$ and $\mathcal R_{>K}(\widetilde G,\widetilde f_0;D)$
are given by Proposition~\ref{prop:mixing} (with $\Omega=D(\varepsilon)$, $G\mapsto\widetilde G$, and $f_0\mapsto\widetilde f_0$);
\item the normal-form remainder error may be taken as
\[
\mathcal E_{\mathrm{nf}}(t,\varepsilon)
\le
C_{\rm err}\,\|\widetilde G\|_{C^1}\,(1+\abs{t}+t^2)\,\|r_\varepsilon\|_\infty;
\]
\item the equilibrium-change error
\begin{align*}
\begin{aligned}
& \mathcal E_{\mathrm{eq}}(\varepsilon):=\\
& \abs{
\int_{D(\varepsilon)\times\T^n}\avgtheta{\widetilde G}(I)\,\widetilde f_0(\theta,I)\,\dd\theta\,\dd I
-
\int_{D(\varepsilon)\times\T^n}\avgtheta{G}(I)\,f_0^D(\theta,I)\,\dd\theta\,\dd I
}
\end{aligned}
\end{align*}
measures the discrepancy between the \emph{angle-averaged equilibrium value in normal-form coordinates}
and the \emph{angle-averaged equilibrium value in the original coordinates}. If $\Phi_\varepsilon^{\rm nf}$ is $C^1$-close to the identity
(as is typically ensured by normal-form lemmas), then $\mathcal E_{\rm eq}(\varepsilon)$ can be further bounded by
$O(\|\Phi_\varepsilon^{\rm nf}-\mathrm{Id}\|_{C^1})$.
\end{itemize}

In particular, if
\[
\|r_\varepsilon\|_\infty\le C_{\mathrm{nf}}\varepsilon e^{-c_{\mathrm{nf}}K(\varepsilon)}
\quad\text{and}\quad
\abs{t}\le \exp(\sigma K(\varepsilon)),\ \ 0<\sigma<c_{\mathrm{nf}}/2,
\]
then
\[
\mathcal E_{\mathrm{nf}}(t,\varepsilon)
\le
C\,\|\widetilde G\|_{C^1}\,\varepsilon
\exp\!\big(-(c_{\mathrm{nf}}-2\sigma)K(\varepsilon)\big),
\]
so this term remains exponentially small in $K(\varepsilon)$ over an exponentially long time window.
\end{theorem}

\begin{proof}
Decompose the initial density as
\[
f_0=f_0^D+f_0^N,\qquad f_0^D:=f_0\1_{D(\varepsilon)},\quad f_0^N:=f_0\1_{N(\varepsilon)}.
\]
Accordingly,
$\langle G\rangle_t=\langle G\rangle_t^D+\langle G\rangle_t^N$ and
$\langle \avgtheta{G}\rangle_0=\langle \avgtheta{G}\rangle_0^D+\langle \avgtheta{G}\rangle_0^N$.

\textbf{(1) A crude bound on the resonant region $N(\varepsilon)$:}
\[
\abs{\langle G\rangle_t^N-\langle \avgtheta{G}\rangle_0^N}
\le
\abs{\langle G\rangle_t^N}+\abs{\langle \avgtheta{G}\rangle_0^N}
\le
2\norm{G}_\infty\int_{N(\varepsilon)\times\T^n}f_0
=
2\norm{G}_\infty\,\mathsf P_{\mathrm{res}}(\varepsilon).
\]

\textbf{(2) Pullback to normal-form coordinates on $D(\varepsilon)$:}
By Lemma~\ref{lem:nf-pullback},
\[
\langle G\rangle_t^D
=
\int \widetilde G(\widetilde\Phi_\varepsilon^t(z))\,\widetilde f_0(z)\,\dd z.
\]

\textbf{(3) Comparing the normal-form flow with the integrable principal flow via the remainder $r_\varepsilon$:}
By Lemma~\ref{lem:nf-error},
\[
\abs{
\int \widetilde G(\widetilde\Phi_\varepsilon^t)\,\widetilde f_0
-
\int \widetilde G(\Psi_\varepsilon^t)\,\widetilde f_0
}
\le \mathcal E_{\mathrm{nf}}(t,\varepsilon).
\]

\textbf{(4) Applying Proposition~\ref{prop:mixing} to the integrable principal flow:}
For $\Psi_\varepsilon^t(\theta,I)=(\theta+t\omega_\varepsilon(I),I)$,
Proposition~\ref{prop:mixing} (on $\Omega=D(\varepsilon)$ with $G\mapsto\widetilde G$ and $f_0\mapsto\widetilde f_0$) yields
\begin{align*}
  \begin{aligned}
&\abs{
\int \widetilde G(\Psi_\varepsilon^t)\,\widetilde f_0
-
\int \avgtheta{\widetilde G}(I)\,\widetilde f_0(\theta,I)\,\dd\theta\,\dd I
}\\
&\le
\frac{C_{\widetilde G,\widetilde f_0}(K(\varepsilon);D(\varepsilon))}{\abs{t}}
+
\mathcal R_{>K(\varepsilon)}(\widetilde G,\widetilde f_0;D(\varepsilon)).
  \end{aligned}
\end{align*}
Finally, we use $\mathcal E_{\mathrm{eq}}(\varepsilon)$ to relate
$\int \avgtheta{\widetilde G}\,\widetilde f_0$ to $\langle \avgtheta{G}\rangle_0^D$.
Combining (1)--(4) gives the result.
\end{proof}

\begin{corollary}\label{cor:nek}
Under the assumptions of Theorem~\ref{thm:main}, suppose we choose $K(\varepsilon)=\lfloor \varepsilon^{-a}\rfloor$ with $a>0$
and that $\|r_\varepsilon\|_\infty\lesssim \varepsilon e^{-c\varepsilon^{-a}}$.
Then for any $0<\sigma<c/2$, whenever $\abs{t}\le \exp(\sigma\varepsilon^{-a})$, one has
\begin{align*}
  \begin{aligned}
\abs{\langle G\rangle_t-\langle \avgtheta{G}\rangle_0}
\le&
2\norm{G}_\infty\,\mathsf P_{\mathrm{res}}(\varepsilon)
+\frac{C_{\widetilde G,\widetilde f_0}(\varepsilon^{-a};D(\varepsilon))}{\abs{t}}\\
&+\mathcal R_{>\varepsilon^{-a}}(\widetilde G,\widetilde f_0;D(\varepsilon))
+
C\norm{\widetilde G}_{C^1}\varepsilon e^{-(c-2\sigma)\varepsilon^{-a}}
+\mathcal E_{\mathrm{eq}}(\varepsilon).
  \end{aligned}
\end{align*}
\end{corollary}

\begin{remark}
 If $G$ and $f_0$ are analytic in $\theta$ (with strip width $\sigma_0>0$), then their Fourier coefficients decay exponentially, and
\[
\mathcal R_{>K}(\widetilde G,\widetilde f_0;D)\lesssim e^{-\sigma_0 K},
\]
so the high-frequency tail can be written explicitly as $e^{-\sigma_0K(\varepsilon)}$.

\end{remark}
\begin{remark}
  A traceable expression for $C_{\widetilde G,\widetilde f_0}(K;D)$ is provided by Proposition~\ref{prop:mixing};
the key quantity is
\[
\gamma_k(D)=\inf_{I\in D}\abs{\grad\big(k\cdot\omega_\varepsilon(I)\big)}.
\]
If $D\omega_\varepsilon$ is uniformly invertible on $D$, then $\gamma_k(D)\gtrsim \abs{k}$,
and consequently the growth of $C_{\widetilde G,\widetilde f_0}$ can be further controlled.
\end{remark}

\section{Conclusion}\label{sec:conclusion}
This paper provides a ``ready-to-use'' quantitative phase-mixing estimate for integrable action--angle flows, with an explicit decay rate of $1/|t|$ and a \emph{traceable} mixing constant: it depends on nonstationary-phase quantities (e.g., the gradient/Hessian of $k\cdot\omega$) and on $W^{1,1}$-norms of the relevant Fourier modes. By combining a resonant/nonresonant decomposition with an analytic normal-form reduction on $D(\varepsilon)$, we embed this mixing estimate into the nearly integrable framework, thereby obtaining an explicit deviation bound for ensemble expectations that holds for arbitrary initial data. The final bound transparently separates the contributions of different error sources: the resonant neighborhood enters through a probability-mass term, while the nonresonant region contributes via the mixing term, the high-frequency Fourier tail, and explicit normal-form remainder and coordinate-change error terms.

A natural continuation of this work is to make the decomposition-based deviation bound fully explicit in $\varepsilon$.
In particular, one may seek geometric estimates for the resonant neighborhood $N(\varepsilon)$ that translate regularity assumptions on the initial marginal density $\rho_0$ and nondegeneracy properties of the frequency map $\omega$ into computable upper bounds for the resonant probability mass $\mathsf P_{\mathrm{res}}(\varepsilon)$, and then optimize the choice of $(K(\varepsilon),\alpha(\varepsilon))$ to balance the resonant-mass term, the $1/|t|$ mixing term, the Fourier tail, and the normal-form remainder error.
A second direction is to refine the present crude treatment of $N(\varepsilon)$ by incorporating resonant normal forms, with the goal of capturing the contribution of slow resonant transport (including Arnold diffusion mechanisms) to ensemble observables.
Finally, it would be interesting to extend the framework beyond the analytic setting (e.g.\ to Gevrey or finitely differentiable Hamiltonians) and to address higher-order statistics, such as quantitative bounds for time correlations and fluctuation results for deterministic nearly integrable flows.
\section*{Declaration of competing interest}
The authors declare that they have no known competing financial interests
or personal relationships that could have appeared to influence the work
reported in this paper.

\section*{Acknowledgements}
The author Yong Li was supported by National Natural Science Foundation of China
(12071175, 12471183 and 12531009).



\end{document}